\documentclass{article} %
\usepackage{iclr2027_conference,times}

\usepackage{amsmath,amsfonts,bm}

\def\eqref#1{equation~\ref{#1}}

\def\1{\bm{1}}

\DeclareMathAlphabet{\mathsfit}{\encodingdefault}{\sfdefault}{m}{sl}
\SetMathAlphabet{\mathsfit}{bold}{\encodingdefault}{\sfdefault}{bx}{n}

\makeatletter
\renewcommand{\eqref}[1]{\textup{\tagform@{\ref{#1}}}}

\usepackage{dashrule}
\usepackage{xcolor}
\usepackage[T1]{fontenc}

\usepackage{tikz}

\usepackage{amsthm}
\theoremstyle{definition}
\newtheorem{definition}{Definition}[section]

\theoremstyle{plain}
\newtheorem{thm}[definition]{Theorem}

\newtheorem{prop}[definition]{Proposition}

\usepackage{booktabs, makecell, multirow}
\usepackage{tikz}
\usetikzlibrary{arrows.meta,positioning}

\definecolor{policyblue}{RGB}{0,114,178}
\definecolor{informedorange}{RGB}{230,159,0}
\definecolor{annealgreen}{RGB}{0,158,115}
\definecolor{greygreen}{RGB}{110,110,110}

\usepackage{algorithm}
\usepackage{algpseudocode}

\makeatother

\usepackage{hyperref}
\usepackage{url}
\usepackage{amsmath}

\ExplSyntaxOn
\NewDocumentCommand{\signseq}{m}
  {
    \texttt{
      \tl_map_inline:nn {#1}
        { \makebox[0.6em][c]{##1} }
    }
  }
\ExplSyntaxOff

\title{New Records for the Hadamard Maximal Determinant Problem}

\author{
Giorgi Butbaia\textsuperscript{1},
Justin Tan\textsuperscript{2,3},
Pragatheeswaran Vipulanandan\textsuperscript{4},
Xiaoyu Huang\textsuperscript{5},\\
\textbf{
Toby Saunders-A'Court\textsuperscript{1},
Lucas Fagan\textsuperscript{1},
Davide Passaro\textsuperscript{1},
Michele Tarquini\textsuperscript{1}
\& Sergei Gukov\textsuperscript{1}
}
\\[0.6em]
\textsuperscript{1}California Institute of Technology, Pasadena, CA, USA\\
\textsuperscript{2}University of Cambridge, Cambridge, UK\\
\textsuperscript{3}London Institute for Mathematical Sciences, London, UK\\
\textsuperscript{4}University of Miami, Miami, FL, USA\\
\textsuperscript{5}Temple University, Philadelphia, PA, USA\\
[0.4em]
}

\iclrfinalcopy %
\begin{document}

\maketitle

\begin{abstract}
The Hadamard maximal determinant problem seeks an $n\times n$ matrix $X$ with entries in $\{\pm 1\}$, which maximizes the determinant $\vert \det X \vert$ for a given order $n \in \mathbb{N}$. We report matrices attaining new record determinants for orders $n\equiv 3\pmod{4}$ between $103$ and $119$, as well as $n=51$. We provide a proof of optimality within a specific family of circulant--block constructions for $X$. For these orders, we recast the Hadamard problem as a sequence reconstruction problem from a pair of periodic autocorrelation functions, subject to certain arithmetic conditions. We describe several learning--based strategies for reconstruction and provide a comparison with classical simulated annealing--style approaches as a benchmark.

\end{abstract}

\section{Introduction}

The \textit{Hadamard maximal determinant problem}, originally posed by \citep{hadamard93} is to find an $n \times n$ $\{\pm 1\}$--matrix with maximum possible determinant for a given order $n$:
\begin{gather*}
    \mathcal{D}(n) = \max\left\{ |\det{X}|~\vert~ X \in \{-1,+1\}^{n\times n}\right\}\,.
\end{gather*}
Na\"ively, the search space for the problem grows superexponentially as $2^{n^2}$. The problem of finding a general construction for all orders $n$ has remained unsolved for over a hundred years. One reason is that the problem naturally divides based on the residue class of $n$ (mod 4). The most progress has been made for $n \equiv 0$ (mod 4) (Proposition 4, \citep{surveyHadamard}), and the least progress for $n \equiv 3$ (mod 4). Certifiably maximal determinant matrices have been found for all orders up to $n=22$ inclusive \citep{orrick2012maxdet}, with the $n=23$ record conjectured to be maximal. Various constructions have been made in the literature for special infinite families of solutions, and there exist certain ad--hoc examples, but no general construction has been found for any residue class.

Geometrically, the determinant of the matrix $X_n$ gives the signed volume of an $n$--dimensional parallelotope in $\mathbb{R}^n$ spanned by the rows/columns of $X_n$, whose vertices lie on the corners of the hypercube with edge length 2. The volume of $X_n$ is maximized when the defining vectors are orthogonal, yielding an elementary upper bound for all orders as $\vert \det X_n \vert\leq n^{n/2}$. This bound is sharp only when the orthogonality condition is met, which, for $n>2$, is only possible when $n \equiv 0$ (mod 4). For other orders, one must find an optimal relaxation of the orthogonality condition, with different arithmetic `defects' to orthogonality arising for each order $n \not\equiv 0$ (mod 4).

The associated Gram matrix \(G_n=X_nX_n^\top\) encodes the geometry of the row vectors of \(X_n\): its diagonal entries are their squared lengths, while its off-diagonal entries are their pairwise inner products. Choosing a suitable normalization, it is not hard to show that the off-diagonal elements of the Gram matrices corresponding to the maximal determinant matrices assume values $\equiv n$ (mod 4) \citep{brenner72}. For $n \equiv 1$ (mod 4), the natural symmetric ansatz $G_n = (n-t)\mathbf{1}_n + tJ_n$, where $J_n$ is the uniform $n \times n$ matrix of $1$s, yields an upper bound on $\det X_n$ which is sharp for at least finitely many matrices of such orders. A similar ansatz exists for $n\equiv 2$ (mod 4) \citep{surveyHadamard}. This symmetry breaks for $n \equiv 3$ (mod 4) --- neither of the uniform choices at $t = -1, 3$ are satisfactory, and the matrices which yield natural upper bounds arrange themselves into a block--diagonal structure over a uniform $(-1)$--background \citep{Ehlich1964a, Ehlich1964b}. We will provide new record--achieving constructions for maximal determinant matrices in orders $n\equiv 3$ (mod 4) that may be realized as \textit{circulant matrices}. In this case, basic Fourier theory states that the problem of constructing the maximal determinant matrices is equivalent to a sequence reconstruction problem from analytically computed target data.

The problem serves as an excellent benchmark for testing the capabilities of reinforcement-learning approaches in extremely sparse--reward environments. Recent successes, as evidenced by the many applications of RL to combinatorial optimization \citep{kool2019attention,kwon2020pomo,chalumeau2023combinatorial,kim2024srt}, motivate us to apply similar techniques to search for maximal determinant targets improving currently known lower bounds \citep{surveyHadamard, orrick2012maxdet}.

\begin{figure}[t!]
    \centering
    \includegraphics[width=0.95\textwidth]{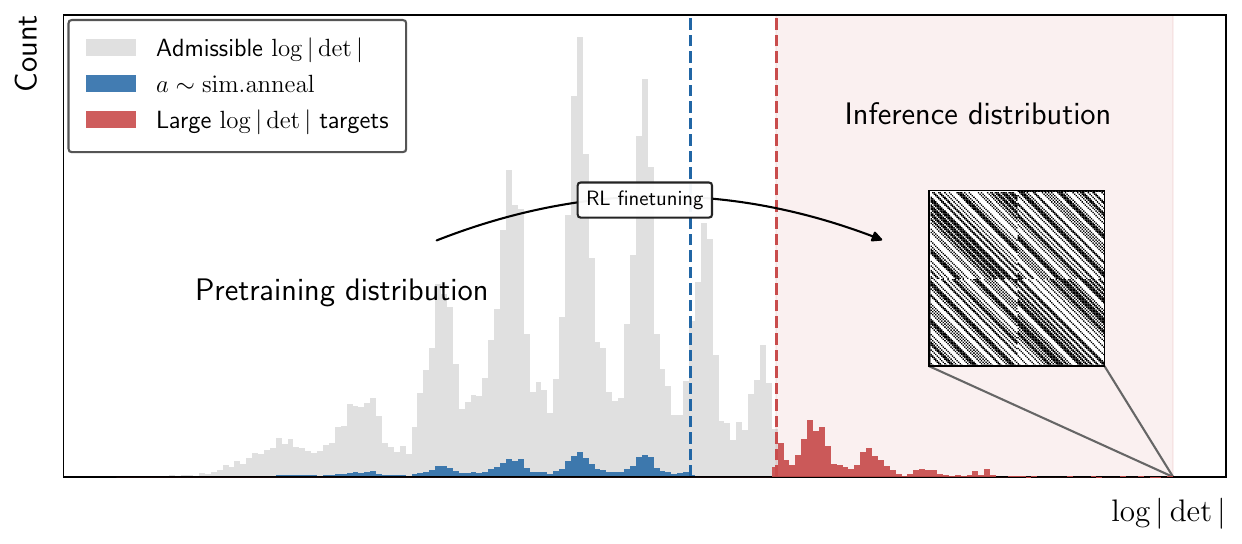}
    \caption{Determinant--maximizing search reformulated as a sequence reconstruction problem using RL. The shaded gray histogram shows the distribution of periodic autocorrelation functions \eqref{eq:paf} in $n=103$, whose determinant may be exactly computed \eqref{eq:determinant}. We train a policy to reconstruct the circulant components of a specific matrix ansatz \eqref{eq:ansatz} on the blue distribution of low--determinant targets, and observe successful generalization to the red distribution of high--determinant targets. Notably, we do not directly supply information about the determinant to the agent.}    
    \label{fig:pafRL}
\end{figure}

We discuss two different reinforcement-learning-based strategies for constructing sequences defining the circulant blocks of the matrices which improve the previously known records. The models are first pretrained on a small set of PAF-sum targets (generated via simulated annealing) and then fine-tuned to generalize to the target distribution of large-determinant matrices (See Fig.~\ref{fig:pafRL}). We demonstrate that this training process successfully transfers to previously unseen PAF-sum targets, including the maximal determinant targets within the families of matrices considered.

\section{Periodic autocorrelations}
Here we reformulate the maximization of the determinant as a sequence reconstruction problem \citep{Kotsireas2019PAF}. In this work, we focus on $\{\pm1\}$-valued $n\times n$ matrices $X$ constructed via $m\times m$ circulant matrices $A$ and $B$, where $m=(n-1)/2$. We consider two main constructions, {\it uncorrelated} $X(A,B)$ and {\it correlated} $X^{\mathrm{corr.}}(A,B)$, defined as follows:
\begin{gather}\label{eq:ansatz}
    X(A,B) := \begin{pmatrix}
    1 & j^T & -j^T \\
    j & A & B \\
    -j & B^T & -A^T
    \end{pmatrix},\quad
    X^{\mathrm{corr}.}(A,B) := \begin{pmatrix}
    1 & j^T & -j^T \\
    j & A & B \\
    j & -B & -A^T
    \end{pmatrix}
\end{gather}
where $j$ is an all-ones vector of size $m$. For any fixed order $n\equiv 3\pmod 4$, we will therefore focus on finding the circulant matrices $A$ and $B$ which maximize the determinant $|\det X(A,B)|$ (or $|\det X^{\mathrm{corr}.}(A,B)|$, whichever is larger). We shall focus on the uncorrelated case here and leave the discussion for the correlated case to Appendix \ref{app:corr}.

An $m\times m$ circulant matrix $A$ is uniquely described by a sequence $\{a_{j}\}_{j=0}^{m-1}$ via $A_{ij} = a_{j-i\pmod m}$. For a circulant matrix $A$ constructed in such way from a sequence $a$, we shall use notation $A = \mathrm{circ}(a)$.

{\bf Uncorrelated $X$}: Let $\mathcal{G}_n = \left\{ XX^T ~\vert~ X \text{ is maximal, } n \equiv 3 \pmod 4\right\}$ denote the set of Gram matrices corresponding to maxdet matrices in orders $n \equiv 3 \pmod 4$. Ehlich's theorem \citep{Ehlich1964b} states the upper bound on the determinants of $G \in \mathcal{G}_n$ is attained by block--structured matrices of the form $\tilde{G}_n = (n-3)\mathbf{1}_n - J + 4\bigoplus_k J_k$, where $J_k$ is the all--ones $k \times k$ matrix. It is not known when this bound is sharp \citep{brent2011maximal, orrick2005maximal}, but it is natural to conjecture that the Gram matrices of high--determinant matrices should satisfy $G_{ij} \in \{-1,+3\}$ for all $i\neq j$. Let $X = X(A,B)$, then the Gram matrix $G$ of $X$ is given by:
\begin{gather}
    G = XX^T = \begin{pmatrix}
        n  & \gamma_0^T & \gamma_1^T \\
        \gamma_0 & C & -j\otimes j \\
        \gamma_1 & -j\otimes j & C
    \end{pmatrix}\,,
\end{gather}
where $\gamma_0 = j + (A-B)j$, $\gamma_1 = -j + (B^T + A^T)j$, and $C = AA^T + BB^T + j\otimes j$. Since $A$ and $B$ are circulant, they may be described using sequences $\{a_j\}_{j=0}^{m-1}$ and $\{b_j\}_{j=0}^{m-1}$, respectively. This allows us to identify the border terms as $\gamma_0 = (1+|a|-|b|)j$, $\gamma_1 = (-1 + |a| + |b|)j$, where $|a| := \sum_{j}a_j$ and $|b| := \sum_{j}b_j$ are segment sums. Furthermore, the diagonal component $C$ is also a circulant matrix:
\begin{gather}
    C = \mathrm{circ}(\{\mathrm{PAF}_a(k) + \mathrm{PAF}_b(k)\}_{k=0}^{m-1}) + j\otimes j
\end{gather}
where $\mathrm{PAF}_a$ and $\mathrm{PAF}_b$ are the periodic autocorrelation functions of sequences $a$ and $b$.

\begin{definition} Let $\{a_j\}_{j=0}^{m-1}$ be a sequence of integers. The {\it periodic autocorrelation function} $\mathrm{PAF}_a(s)$ is defined as:
\begin{gather}\label{eq:paf}
    \mathrm{PAF}_a(s) = \sum_{j=0}^{m-1}a_{j}a_{j+s\pmod m}\, \quad s \in \{0,\ldots,m-1\}~.
\end{gather}
\end{definition}
The Ehlich Gram constraint $G_{ij}\in\{-1,+3\}$ for all $i\neq j$ imposes the following condition on the sum of the periodic autocorrelation functions:
\begin{gather}\label{eq:paf_constraint}
    \forall s>0\colon p(s):= \mathrm{PAF}_a(s) + \mathrm{PAF}_b(s)\in\{-2,+2\}\,,
\end{gather}
where $p(0) = 2m = n-1$. We will refer to such functions as {\it admissible}. By summing the function $p(s)$ over $0\leq s<m$, the Ehlich constraint translate to a Diophantine constraint on the segment sums \citep{djokovic2011new, williamson2026generating}:

\begin{prop} Let $N$ be the number of times $p(s) = +2$ for $s\neq 0$, then:
\begin{gather}\label{eq:diophantinePAF}
    |a|^2 + |b|^2 = 4N+2\,.
\end{gather}
\end{prop}
\begin{proof} The proof of the original statement can be found in \citep{djokovic2011new}. We provide a proof for our parametrization in Appendix~\ref{app:paf}.
\end{proof}

Recall our goal is to reconstruct high--determinant matrices from PAF targets. First we identify pairs of segment sums $(|a|,|b|)$ satisfying \eqref{eq:diophantinePAF}. Next, we need to be able to compute the determinant associated to an admissible joint PAF target $p(s)$ and enumerate the ones exceeding the current record \citep{orrick2012maxdet}. We proceed by decomposing $X$ into its Fourier modes. In particular, let $A$ be a circulant matrix, then its normalized eigenvectors are given by:
\begin{gather}
    v_t = \frac{1}{\sqrt{m}}(1,~\zeta_m^t,~\dots, \zeta_m^{(m-1)t})^T,\quad 0\leq t<m\,,
\end{gather}
where $\zeta_m := \exp(2\pi i/m)$ is the primitive $m$-th root of unity. It follows that the eigenvalues are $\lambda_t = \sum_{k=0}^{m-1}a_k \zeta_m^{kt}$. Let $F = (v_0,\dots, v_{m-1})$, then $F$ diagonalizes $A$ into Fourier modes:
\begin{gather}
    F^\dagger A F = \mathrm{diag}(\widehat{a}_0,\dots, \widehat{a}_{m-1})\,.
\end{gather}
By taking a unitary transform of $X$ using $U:= 1 \oplus F \oplus F$, we obtain a decomposition:
\begin{gather}
    U^\dagger XU \cong \widehat{X}_0 \oplus \bigoplus_{t>0}^{m-1} \widehat{X}_t,\quad\text{(up to permutations)}
\end{gather}
where we reorder the basis to group the $t=0$ modes together, and the respective $t>0$ modes for $A$ and $B$ together. Note that $j$ only couples to the zero-modes, so $\widehat{X}_0$ is the $3\times 3$ matrix corresponding to the zero-mode of $X$ incorporating the border terms:
\begin{gather}
    \widehat{X}_0 = \begin{pmatrix}
    1 & \sqrt{m} & -\sqrt{m} \\
    \sqrt{m} & \widehat{a}_0 & \widehat{b}_0 \\
    -\sqrt{m} & \widehat{b}_0 & -\widehat{a}_0
    \end{pmatrix}\,,
\end{gather}
whereas the higher modes $\widehat{X}_{t>0}$ are given by:
\begin{gather}
    \widehat{X}_{t>0} = \begin{pmatrix}
    \widehat{a}_t & \widehat{b}_t \\
    \overline{\widehat{b}}_t & -\overline{\widehat{a}}_t
    \end{pmatrix}\,.
\end{gather}

Since $U$ is unitary, it preserves determinant, therefore, we obtain:
\begin{gather}\label{eq:determinant}
    \det{X} = \det{X_0}\prod_{t>0}^{m-1}\left(|\widehat{a}_t|^2+|\widehat{b}_t|^2\right)\,.
\end{gather}
The zero-mode determinant is straightforward (it depends only on zero modes of $A$ and $B$), however, the higher modes may be related to the Fourier transforms of the respective PAFs in the following manner: using the Wiener-Khinchin theorem \citep{wiener1930generalized, khintchine1934korrelationstheorie}, the power spectrum of a sequence $(a_0,\dots, a_{m-1})$ is the DFT of its periodic autocorrelation function:
\begin{gather}
    |\widehat{a}_t|^2 = \sum_{s=0}^{m-1}\mathrm{PAF}_a(s)\zeta_m^{st}\,,
\end{gather}
therefore, $|\widehat{a}_t|^2+|\widehat{b}_t|^2$ may be rewritten purely in terms of the Fourier transform of the PAF-sum function $p(s)$ as: $\widehat{p}_t = |\widehat{a}_t|^2+|\widehat{b}_t|^2$. The Fourier transform $\widehat{p}_t$ is given by:
\begin{gather}
    \widehat{p}_t = 2m + 2\sum_{s=1}^{(m-1)/2}p(s)\cos\left(\frac{2\pi s t}{m}\right)\,,
\end{gather}
therefore, in order to evaluate the determinant of the matrix yielding a proposed joint PAF target, we may do so directly using $\hat{p}_t$, without requiring its reconstruction \eqref{eq:paf_constraint} into sequences $a$ and $b$.

\begin{figure}[hb!]
    \centering
    \includegraphics[width=1.0\textwidth, trim = {0.0cm 0.5cm 0.0cm 0.1cm}]{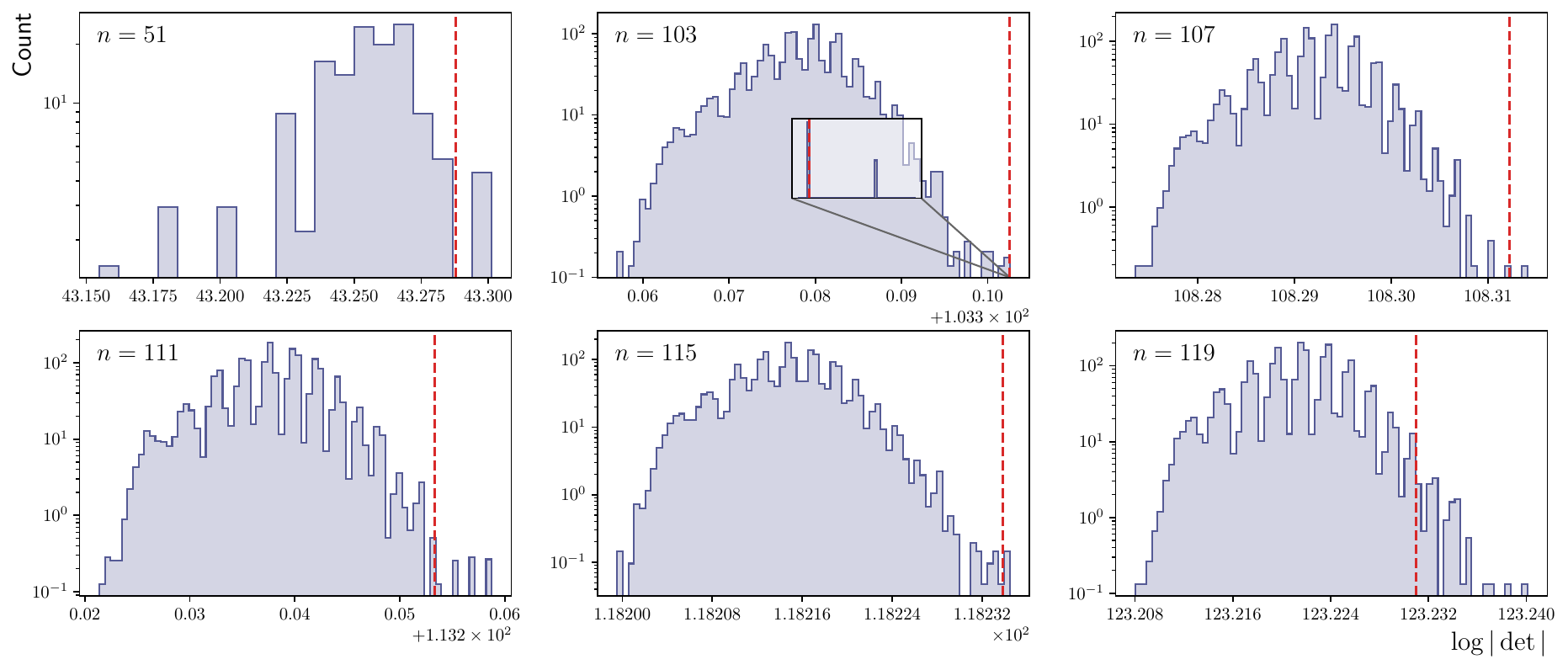}
    \caption{Distribution of $\log|\det|$ values of all admissible PAF-sum targets computed using \eqref{eq:determinant}. The dashed red line (\protect\textcolor{red}{\hdashrule[0.5ex]{2em}{1pt}{2pt 2pt}}) denotes the current record in \citep{orrick2012maxdet}.}%
    \label{fig:paftargets}
\end{figure}

Once we have a list of candidate PAF-sums $\{p_i(s)\}_i$ achieving maximal possible determinant in the uncorrelated ansatz class $X = X(A,B)$ \eqref{eq:ansatz}, we want to reconstruct the original sequence pair $(a,b)$. This is an underdetermined system of quadratic equations in $2m$ variables on the product of boolean cubes $\{\pm 1\}^{m}\times \{\pm 1\}^m$. Equivalently, we would like to minimize the following objective:
\begin{gather}\label{eq:energy}
    \min_{a,b \in \{\pm 1\}^m} H(a,b) = \min_{a,b \in \{\pm 1\}^m} \sum_{s=1}^{(m-1)/2}\bigg|\mathrm{PAF}_a(s) + \mathrm{PAF}_b(s) - p(s)\bigg|^2\,,
\end{gather}
where $s=0$ is omitted as it gives a constant contribution and summation indices end at $(m-1)/2$ as $p(s)$ is reflection--symmetric about the midpoint.

In order to construct a distribution of all PAF-sum targets w.r.t. their expected determinants, we first exhaustively search through the solutions of the Diophantine equation \eqref{eq:diophantinePAF} and for each candidate segment sum pair $(|a|,|b|)$, we sweep through all possible PAF-sum targets. Given that a specified segment sum pair specifies the number of spots where $p(s)=2$ at $s\neq 0$ and noting that $p(s) = p(-s)$, it is immediately obvious that the total number of all possible PAF targets is:
\begin{gather}
    \#\,\left\{\text{admissible targets with segment sum pair}~(|a|,|b|)\right\} = \begin{pmatrix}
    \displaystyle\frac{m-1}{2} \\
    \displaystyle \frac{|a|^2+|b|^2-2}{8}
    \end{pmatrix}\,,
\end{gather}
where the dimension which has the largest number of admissible targets is $n = 119$, which has $475020$ admissible targets, which is small enough to be enumerated completely. Thus reconstructing a maximal PAF-sum target immediately implies a verifiable certificate of maximality within the family $X(A,B)$ with $p(s\neq 0)\in\{-2,2\}$. The distribution of all such targets is shown in Fig.~\ref{fig:paftargets}. We find empirically that the highest--determinant candidates degenerate to a single pair of segment sums at all orders. The correlated ansatz $X^{\mathrm{corr.}}$ has a similar construction, described in Appendix~\ref{app:corr}.

\section{Maximal determinants as sequence reconstruction}\label{sec:seq_recon}
Finding the maximal determinant representatives of the families \eqref{eq:ansatz} may be regarded as a reconstruction problem which attempts to construct a pair of sequences  $(\{a_j\}_{j=0}^{m-1}, \{b_j\}_{j=0}^{m-1})$ satisfying the relevant constraints depending on whether they parametrize $X(A,B)$ or $X^{\mathrm{corr.}}(A,B)$. While $X(A,B)$ matrices allow every pair of sequences with values in $\{\pm 1\}$ minimizing the energy objective \eqref{eq:energy}, $X^{\mathrm{corr.}}(A,B)$ additionally requires two more conditions (See appendix~\ref{app:corr}):
\begin{gather}
    |a||b| = m,\quad a*c = 2,\quad\text{where}~c_t = b_t+b_{-t}\,.
\end{gather}
In both cases, we define the state space $\mathcal{S}(|a|,|b|,m)$ as the set of all pairs of sequences $(a,b)$ of length $m$ with segment sums $\sum_j a_j = |a|$ and $\sum_j b_j = |b|$. The initial state $s_0$ is sampled uniformly from the state space. Since we are interested in reconstructing PAF-sum targets using sequences with pre-determined segment sums $|a|$ and $|b|$, we define the action space $\mathcal{A}(|a|,|b|,m)$ to be the set of all moves which swap two distinct-valued entries within sequence $a$ or sequence $b$. From this, it immediately follows that the state and action space sizes scale as:
\begin{gather}\label{eq:stateSize}
    \left|\mathcal{S}(|a|,|b|,m)\right| = \begin{pmatrix}m \\ \displaystyle \frac{m+|a|}{2}\end{pmatrix}\begin{pmatrix}m \\ \displaystyle \frac{m+|b|}{2}\end{pmatrix},\quad |\mathcal{A}(|a|,|b|,m)| =  \frac{2m^2-|a|^2-|b|^2}{4}\,.
\end{gather}
Typical scales of state spaces that we consider are shown in Table~\ref{tab:stateSpaceScaling}. Due to the extreme combinatorial explosion, the performance of our methods depends significantly on the design choices for learning useful representations and the reward signals.

\subsection{Arithmetic and hierarchical reconstructions}
At the orders we are interested in, as indicated in Table~\ref{tab:stateSpaceScaling}, the search space is prohibitively large for any direct reconstruction algorithm relying on either pure RL or classical techniques. We therefore break down the reconstruction problem into two steps:
\begin{enumerate}
    \item Given a PAF--sum target $p(s)$, find the space of all feasible sequences $\{a_j\}_{j=0}^{m-1}$ such that the residual:
    \begin{gather}\label{eq:pafResidual}
        \mathrm{PAF}_b(s) = r(s) := p(s) - \mathrm{PAF}_a(s)
    \end{gather}
    is algebraically unobstructed, where an obstruction is defined below.
    \item Given an algebraically unobstructed residual $r(s)$, attempt to perform an exhaustive search for a sequence $\{b_j\}_{j=0}^{m-1}$ solving \eqref{eq:pafResidual}.
\end{enumerate}
Motivated by the use of cyclotomic fields in the recent works \citep{djokovic2011new, kotsireas2025quaternary}, we also consider a projection of a sequence $\{a_j\}_{j=0}^{m-1}$ to a cyclotomic field $K_d := \mathbb{Q}(\zeta_d)$ for $d\mid m$, and rewrite the residual equation \eqref{eq:pafResidual} as a norm equation:
\begin{gather}\label{eq:pafResidualCyclotomic}
    b_db_d^* = r_d\,.
\end{gather}
From this identification, it follows that all obstructions to solvability of the field-norm equation may be detected by studying the class of $r_d$ in Galois cohomology $H^2(K_d/K_d^+)$ (See Appendix~\ref{app:arithmetic} for more details).
\subsection{Reconstructions through lifts}\label{ssec:oracle}
Let $m = (n-1)/2$ be a non--prime odd number and $d\mid m$ be one of the factors such that $d>2$ and $d\neq m$. While the arithmetic handles obstructions at the level of $\mathbb{Q}(\zeta_d)$, where $\zeta_d$ is a $d$--th root of unity, we may also consider projections of the form:
\begin{gather}
    \pi_d\colon \mathbb{Z}[C_m]\rightarrow \mathbb{Z}[C_d]\,,
\end{gather}
here $C_k$ is the cyclic group of order $k$. We will first attempt to construct a set of potential solutions in $\mathbb{Z}[C_d]$ and then lift them to the full solution in $\mathbb{Z}[C_m]$. More specifically we run two exhaustive searches:
\begin{itemize}
    \item The coefficients of $\pi_d(b)$ have values in $\{-k,-k+2,\dots, k\}$ where $k = m/d$. Projections preserve segment sums and map PAFs to PAFs. Therefore, we exhaustively construct a set:
    \begin{gather}\label{eq:proj}
        \mathcal{P}_d(r, |b|) = \left\{g\in \{-k,-k+2,\dots, k\}~\vert~ |g| = |b|,~\mathrm{PAF}_g = \pi_d(r)\right\}\,.
    \end{gather}
    \item If $\mathcal{P}_d(r,|b|) \neq \emptyset$, we search for all possible $\{\pm 1\}$ sequences with prescribed fiberwise sums from $g\in \mathcal{P}_d(r,|b|)$:
    \begin{gather}\label{eq:lift}
        \mathcal{B}_d(g) = \left\{b\in\{\pm 1\}^m~\bigg\vert~ \forall u\in\mathbb{Z}_d\colon \sum_{j=0}^{k-1}b_{u+jd} = g_u\right\}\,.
    \end{gather}
    The elements of $\mathcal{B}_d(g)$ are the candidate reconstructions of the PAF equation. We then exhaustively check $\mathrm{PAF}_b = r$ for every $b\in \mathcal{B}_d(g)$.
\end{itemize}
We massively parallelize the checks \eqref{eq:proj} and \eqref{eq:lift} over our custom CUDA kernels. This algorithm defines an \textit{oracle} $\mathsf{B}$ that given a candidate sequence $a \in \{\pm 1\}^m$ and a target joint autocorrelation function $p(t) \in \{\pm 2 \}^m, \, t > 0$, is able to provide a valid partner $b$ satisfying the conjectured maximal--determinant condition \eqref{eq:paf_constraint}, or certify that no such partner exists, i.e. 
\begin{equation}\label{eq:oracle}
    \mathsf{B}(p;A) := \mathbf{1}\left[\exists \, b \in \{\pm 1\}^m : \mathrm{PAF}_A + \mathrm{PAF}_B = p \right]~.
\end{equation}
The problem then reduces to construction of $A$, and the effective horizon length halves. %

\subsection{Arithmetic RL}\label{sec:arithmeticRL}
The major difficulty of two-stage reconstruction is finding a large number of sequences $\{a_j\}_{j=0}^{m-1}$ which pass the field norm obstructions \eqref{eq:cohomology}. Although passing these obstructions provides a relatively sparse signal, the decomposition of the global obstruction problem to local tests over relevant places $0\neq v\in\mathrm{Spec}~\mathcal{O}_{K_d^+}$ allows us to construct a potential which provides useful feedback for guiding an agent. Therefore, we build an agent, trained using an algorithm similar to that of our RL setup, described in Sec.~\ref{sec:jointRL}, whose main objective is to sample a large number of unique sequences with unobstructed residuals for any prescribed PAF-sum target. In particular, we first pretrain the model on a dataset of PAF-sum targets produced using simulated annealing via algorithms ~\ref{alg:samplePretrainingData}, and \ref{alg:pretrainAgent}. We then finetune the model using curriculum scheduled reinforcement learning (See Alg.~\ref{alg:curriculumFinetuning}), after which, we perform inference on previously unseen PAF-sum targets. We illustrate this procedure on Fig.~\ref{fig:pafRL}.

\begin{figure}[h!]
    \centering
    \begin{tikzpicture}
        \node[anchor=south west, inner sep=0] (graph) at (0,0) {
            \includegraphics[width=0.85\linewidth]{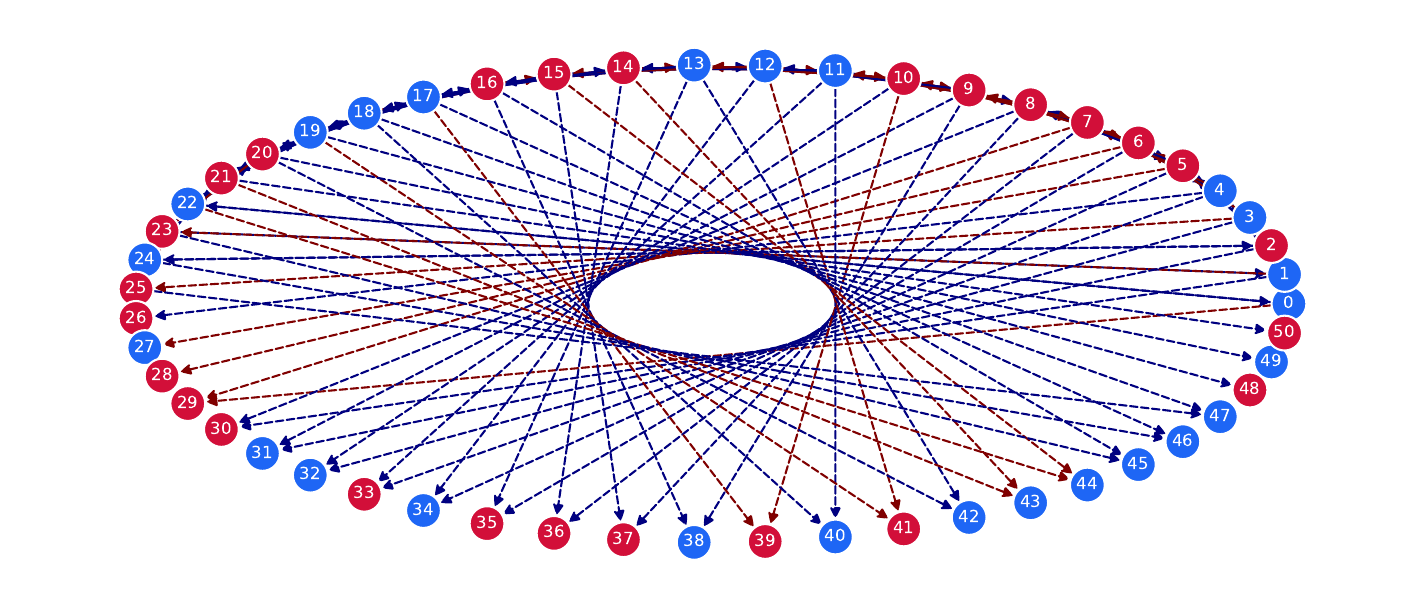}};

        \begin{scope}[
            x={(graph.south east)},
            y={(graph.north west)},
            shift={(graph.south west)}]
            \node[text opacity=1,
                rounded corners=2pt,
                inner sep=5pt] (formula) at (0.97,0.90) {
                $\displaystyle
                \alpha_{ij} = \frac{\exp\left(w_e^\top e_{ij}+b_e\right)}
                {\sum_{v\in\mathcal N(i)}
                 \exp\left(w_e^\top e_{iv}+b_e\right)}$};
            \draw[-{Stealth[length=3mm]},line width=0.9pt,black!90] (formula.south) to[out=-90,in=30] (0.65,0.55);
        \end{scope}
    \end{tikzpicture}
    \caption{Single-PAF GNN model. Nodes represent individual sequence elements $a_j$, where {\color{red}red} nodes denote $a_j=+1$, while {\color{blue}blue} nodes denote $a_j=-1$. The edges are colored based on the values of message weights $\alpha_{ij}$. For each node, only top 2 incident edges with the largest weights are shown.}
\end{figure}

{\bf Single-PAF GNN model:} In order to ensure equivariance, we represent the sequence $\{a_{j}\}_{j=0}^{m-1}$ as a complete graph $\mathcal{K}_m$, where to each edge $e_{ij}$, we assign features:
\begin{gather}
    e_{ij} = \begin{pmatrix}
        p_{d(i,j)} \\
        r_{d(i,j)}
    \end{pmatrix}\,,
\end{gather}
where $d(i,j) = \min(|i-j|, m - |i-j|)$, while the node features are given by $x_i = \phi_{\mathrm{sign}}(a_i) + \phi_{\mathrm{tgt}}(r,t)$, where $\phi_{\mathrm{sign}}$ and $\phi_{\mathrm{tgt}}$ are encoding layers for node and target features, respectively. The edge features, computed via an $\mathrm{MLP}_\theta(h_i,h_j)$, where $h_i,h_j$ are the node features from the final layer, are then used as logits for selecting a pair of nodes to swap.

\begin{figure}
    \centering
    \includegraphics[width=1.0\textwidth, trim = {0cm, 0.8cm, 0cm, 0cm}]{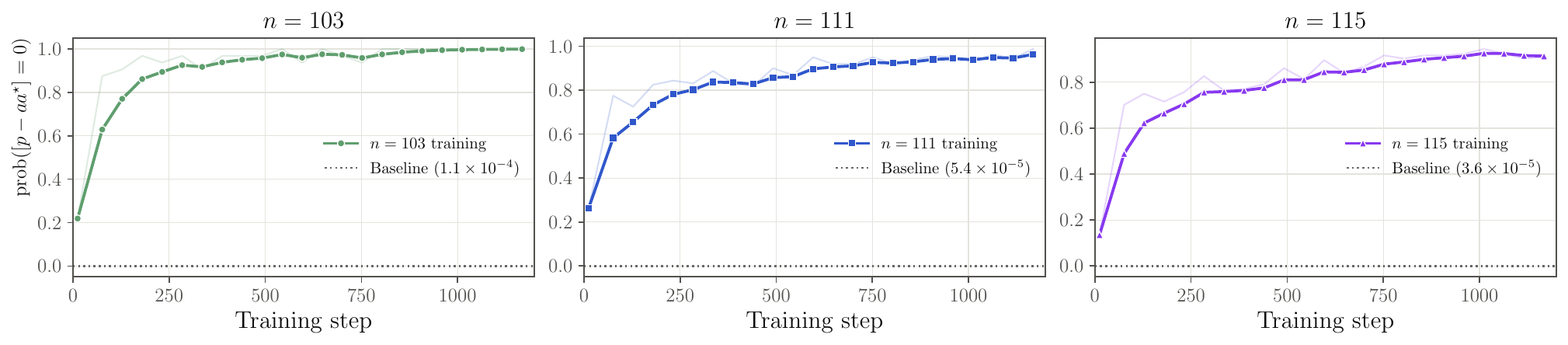}
    \caption{Success rate curves for different orders of the RL agent throughout the training. The success is defined as the environment terminating before truncation due to encountering a state with unobstructed residual relative to maximal determinant PAF-sum target. The dotted horizontal lines denote the baseline probabilities defined as the likelihood of randomly sampled sequence $a$ with $|a| = 1$ having unobstructed residual for the maximal determinant target.}\label{fig:arithmeticTraining}
\end{figure}

{\bf Reward shaping:} We are after two key objectives: given a PAF--sum target $p(s)$:
\begin{itemize}
    \item Find a sequence $a$, such that the Galois cohomology class of the projection of the residual $r_d=p_d - a_da_d^*\in K_d$ vanishes in $H^2(K_d/K_d^+)$ for all $d\mid m$.
    \item Maximize the number of unique sequences $a$ sampled using RL.
\end{itemize}
Based on this, we first define a potential based on the local field norm obstructions from Thm.~\ref{thm:obstructions}:
\begin{gather}
    V_d(a) = -\frac{1}{[K_d^+:\mathbb{Q}]}\sum_{\substack{\text{Embedding}\\
    \sigma\colon K_d^+\rightarrow\mathbb{R}}}\mathbf{1}(\sigma(r_d)<0) - \sum_{\substack{v\in\mathrm{Spec}~\mathcal{O}_{K_d^+}\\v\neq 0}}\mathbf{1}\left((r_d,\Delta_d)_v \neq 1\right)\,,
\end{gather}
which allows us to shape the reward function as:
\begin{gather}
    r(s',a,s) = r_{\mathrm{term}}(s',a,s) + r_{\mathrm{novel}}(s',a,s) + \sum_{d\mid m,~d>2}\left(\gamma V_d(s')-V_d(s)\right)\,,
\end{gather}
where $r_{\mathrm{term}}(s',a,s)$ provides a constant reward of $c_{\mathrm{term}}$ if the residual associated to $s'$ is unobstructed, otherwise it's zero, and $r_{\mathrm{novel}}(s',a,s) = \mathbf{1}\left(s'\notin \mathcal{D}_{\mathrm{seen}}(p)~\wedge~ [p-s's'^*] = 0\right)$ incentivizes previously unobserved sequences $a$ for a given PAF-sum target $p$, where we maintain a set $\mathcal{D}_{\mathrm{seen}}(p)$ of all sequences that have been observed to have unobstructed residuals w.r.t. target $p$. In order to further incentivize the novelty of the constructed sequences, we also add additional reward term coming from random network distillation \citep{burda2019exploration}. The training curves of this model are shown in Fig.~\ref{fig:arithmeticTraining}, while the quantitative reconstruction results are in Fig.~\ref{fig:pass_at_k_simpleRL} and Table~\ref{tab:rl_efficiency}.

\subsection{Arithmetic--free baseline}\label{sec:arithmeticFree}
To understand the utility of the obstruction--theoretic machinery constructed above, we develop a simple baseline that makes minimal assumptions about the underlying structure of the binary sequence $A$. 

\subsubsection{Environment}
The core idea is to work with gaps between consecutive $+1$ elements in the sequence $A$, thereby quotienting out part of the shift symmetry of the autocorrelation function. Denote by $S_A = \left\{ 0 = c_0 < c_1 < \cdots < c_{k_A-1} < m\right\}\subseteq \mathbb{Z}_m$ the support of the $+1$ elements in $A$, and the cardinality $k_A := \vert S_A \vert = (\sigma_A + m)/2$.
A valid trajectory lies in the set
\begin{equation}
    \mathcal{S} := \left\{(g_1, \ldots, g_{k_A-1}) \mid g_{i+1} = (c_{i+1} - c_i) \in \{1, \ldots, m-k_A + 1\} ~\wedge~ \sum_i g_i = m-1\right\}~,
\end{equation}
which is produced autoregressively by a target--conditioned policy $\pi_{\varphi}(\cdot \vert s,p)$ (we will elide $p$ in the sequel) and the state space consists of all valid prefixes. Note $\vert \mathcal{S} \vert = \binom{m-1}{k_A - 1}$ is combinatorially large. The action space consists of a discrete choice of gap $g_i$ at each step along the trajectory. Note the final gap is forced and we rotate the stream such that $c_0 = 0$. This reduces the number of representatives of the shift orbit for each trajectory in $\mathcal{S}$ from $m$ to $k_A$. Suppose we have committed the first $t$ gaps and have a partial reconstruction of $A$, with the last committed position $c_t$. Define $R_t = m - c_t$ as the remaining length to be filled at step $t$ and $L = k_A - t $ as the remaining length of the trajectory. Then the next gap $g_{t+1} \in \{1, \ldots, R_t-L_t+1\}$, and we mask actions accordingly.

Our rewards can only depend on $A$ since we have no access to $B$ during the rollout. Define $\Delta(s) := \mathrm{PAF}_A + \mathrm{PAF}_B - p$. Parseval's theorem says that the Fourier transform pair $(\Delta, \hat{\Delta}_t)$ satisfies
\begin{equation*}
    \sum_s \Delta(s)^2 = \frac{1}{m}\sum_t \vert\hat{\Delta}_t\vert^2~.
\end{equation*}
 A necessary and sufficient condition for $(A,B)$ to be a reconstruction is that $\hat{\Delta} \equiv 0$ for a complete trajectory. From the definition of the Hamiltonian, $\sum_t \vert \hat{\Delta}_t \vert^2  \propto H$. Noting $\hat{\Delta}_t = \vert \hat{a}_t\vert^2 + \vert \hat{b}_t \vert^2 - \hat{p}_t$, we see that non--negativity of $\vert\hat{b}_t\vert^2$ requires the lower bound
\begin{equation}
    \underline{H}(A) := \sum_t \left[ \max(0,\vert \hat{a}_t \vert^2 - \hat{p}_t)\right]^2 \leq \sum_t \vert\hat{\Delta}_t \vert^2~,
\end{equation}
on the Hamiltonian to vanish in order for $A$ to admit a partner. We encourage this behaviour via the following reward, awarded on terminal states $t=T:=k_A-1$,
\begin{equation}
    r(s_t) = \mathbf{1}[t=T] \cdot \chi(A) := \mathbf{1}[t=T] \cdot \mathbf{1}\left[ \underline{H}(A) = 0\right]~.
\end{equation}
While this is a necessary condition for the existence of a partner to $A$, the reward only accounts for feasibility of $A$ without any reference to $B$ (realisability). However, provided the policy is sufficiently diverse, this provides sufficient signal for the trained policy to efficiently sample realisable pairs $(A,B)$. We attempted to augment this with a terminal reward from the oracle \eqref{eq:oracle} on the complete $A$, $r(s_t) \mapsto r(s_t) + \lambda \cdot \mathbf{1}[t=T] \cdot \mathsf{B}(p; A(s_t))~, \, \lambda \in \mathbb{R}^+$, but the oracle rewards appear too sparse to have a significant effect on performance.

\subsubsection{Model}
To motivate our design, note our strategy for reconstructing the pair of sequences which yield a target joint autocorrelation function $p(t)$ is to sample in parallel from the learned policy $\pi_{\varphi}$ to yield $A$, then obtain $B$ or a certificate of infeasibility of $A$ from an oracle, which can be constructed for orders $n$ where $m = (n-1)/2$ is not prime. Note that the question of whether a given $A$ has a partner $B$ yielding a given target is deterministic, and all the randomness arises from the policy. Define the pass rate for the $i$--th target as $\rho_1^{(i)} := \mathbb{P}_{A \sim \pi_{\varphi}}[\mathsf{B}(p^{(i)}; A) = 1]$, then the pass@1 rate is averaged over the validation set $V$ $\rho_1 = \vert V\vert^{-1} \sum_{i \in V} \rho_1^{(i)}$. It will also be useful to record the probability that at least a single feasible example is recovered from $k$ samples --- the pass@$k$ rate $\rho_k$\footnote{Defined similarly per--target as $\rho_k^{(i)} = \mathbb{P}_{A_1,\ldots,A_k \sim \pi_{\varphi}}\left[ \exists j : \mathsf{B}(p; A_j) = 1\right] = 1 - (1-\rho_1^{(i)})^k$.}, of which an unbiased estimator for $n_i \geq k$ samples is:
\begin{equation}\label{eq:unbiased_est_rho_k}
    \hat{\rho}_k = \frac{1}{\vert V\vert} \sum_{i \in V} \left(1 - \binom{n_i - n_i^+}{k} \bigg/ \binom{n_i}{k}\right)
\end{equation}
with $n_i^+$ the number of successful reconstructions for the $i$--th target. A given target $p$ may be realised by a very large set of $A$, and all realisations of $p$ yield the same determinant and are equally good. Because the agent is only rewarded for feasibility, policy diversity is crucial to efficiently certify a given target autocorrelation as attainable --- i.e. the probability mass of the policy should be spread over many different feasible sequences. The oracle \eqref{eq:oracle} itself is deterministic, so sampling identical non--realisable sequences multiple times is detrimental to our goal.

For this reason, we regularize the objective with a KL term between the policy and a uniform prior over legal completions of a partial trajectory, given as a dense reward at each timestep. The trained policy is then defined as the conformal shift of $\pi_0$ (with actions $a$ zero--indexed):
\begin{equation*}
    \pi_0(a_t \vert s_t) := {R_t-a_t-2 \choose L_t-2} \bigg/ {R_t-1 \choose L_t-1}~, \quad \pi_{\varphi}(a_t \vert s_t) \propto \pi_0(a_t \vert s_t) \exp(-f_{\varphi}(s_t;a_t))~.
\end{equation*}
Here the policy head $f_{\varphi}$ is a standard autoregressive transformer module, optimized using PPO \cite{schulman2017proximal} without an entropy bonus. We describe the architecture used for this baseline model in Appendix \ref{app:baseline_arch}. Written over the entire sampled trajectory, the objective we are implicitly maximizing reads as,
\begin{equation}\label{eq:implicit_obj}
    J(\varphi) = \mathbf{E}_{\pi}\left[\sum_t r_t
    \right]-\alpha \cdot \textsf{KL}\left( P_{\pi} \Vert P_0\right) = \mathbf{E}_{\pi}\left[ \sum_t \left(r_t - \alpha \left(\log \pi_{\varphi}(a_t \vert s_t) - \log \pi_0(a_t \vert s_t)\right)\right)\right]~.
\end{equation}
Here $P_{\pi}$ and $P_0$ denote the joint probability of the sampled trajectory under the policy and prior respectively. We found this simple regularization of the return to yield a significantly higher validation pass@$k$ than standard PPO, which tended to have a high pass@1 on training targets but very low pass@$k$--statistics on the validation set --- the actor would place an inordinate amount of probability mass on a relatively small number of feasible sequences $A$ with many possible completions. This is not surprising as maximising the expected return without regularisation towards a prior penalises the policy for placing any probability mass onto non--optimal sequences \citep{korbak-etal-2022-rl}.

\subsection{Results}
To demonstrate the utility of RL in this setting, we sample uniformly permissible sequences $A$ (with $k_A$ positive elements) and run a Monte Carlo estimate of the pass@1 rate for this `random' policy, with a fixed number of samples per target at each order $2^{24}$ at $n=103$ and $2^{25}$ at $n=111,115$, and exhibit the results in Table \ref{tab:rl_efficiency}. We compare this against pass@1 on the validation set for the trained arithmetic--free and arithmetic--aware policies. The trained arithmetic--free baseline policy increases the sampling efficiency by over three orders of magnitude for each $n$, while taking the arithmetic structure of the sequence $A$ into account yields roughly another order of magnitude improvement.

For orders $n=103,111,115$, the policy obtained using both methods is able to realize the sequences $(A,B)$ corresponding to the highest determinant predicted, $\max_p \det X_0 \prod_{t > 0}\hat{p}_t$, via an exhaustive scan of the autocorrelation targets $p(t)$. The results for the record--setting sequences and the efficiency of our scan on the validation set (the highest--$\vert \det X\vert$ 5\% of targets) are reported in Table \ref{tab:record_matrices} in Appendix \ref{app:further_results}. Note that not every target autocorrelation may be guaranteed to admit a realizing pair $(A,B)$, and, while being close to 1, the coverage reported on the validation set should be interpreted as a lower bound on the realizable fraction.%

\begin{figure}
    \centering
    \includegraphics[width=1.0\linewidth, trim={0cm, 0.5cm, 0cm, 0cm}]{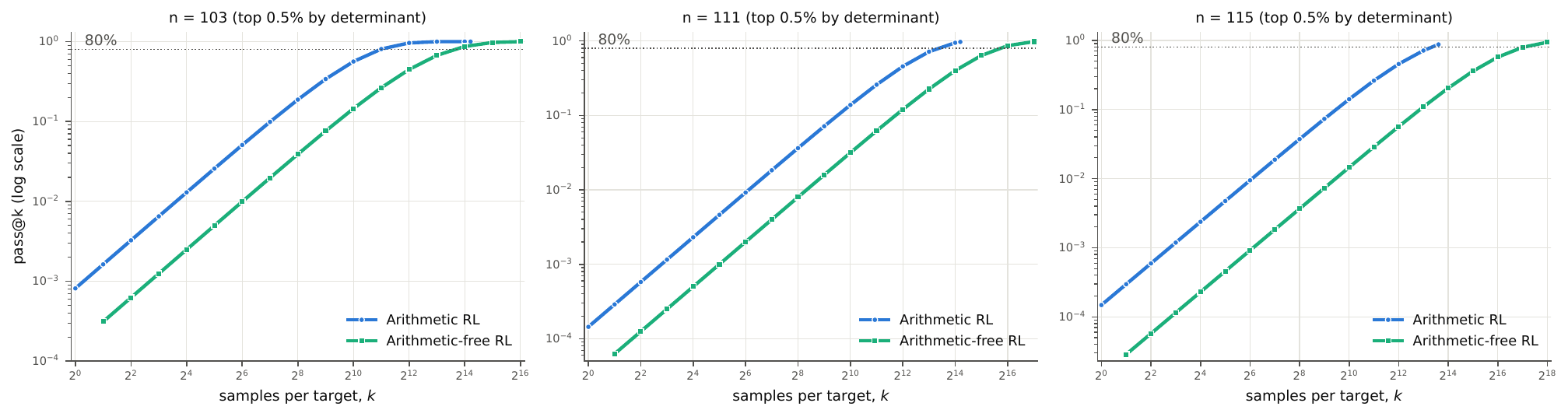}
    \caption{pass@$k$ curves for the final trained policies, computed using the unbiased estimator for $\rho_k$ \eqref{eq:unbiased_est_rho_k} using both Arithmetic-free \ref{sec:arithmeticFree} and Arithmetic RL \ref{sec:arithmeticRL}.}
    \label{fig:pass_at_k_simpleRL}
\end{figure}

From the pass@$k$ rate evaluated for the entire validation set and the top--end, exhibited in Figure \ref{fig:pass_at_k_simpleRL} and Table \ref{tab:validation_samples}, we observe the difficulty of reconstruction is relatively insensitive to the determinant of the predicted target. Further, it appears each increment of $+4$ in the order $n$ roughly corresponds to a two--fold increase in reconstruction difficulty. 

\begin{table}
\caption{Increase in sampling efficiency due to the trained policy. Note the advantage grows with the order.}\vspace{0.1cm}
\centering
\begin{tabular}{rlrr}
\toprule
$n$ & Random policy pass@1 & Arithmetic-free / rand. $(\uparrow)$ & Arithmetic / rand. $(\uparrow)$ \\
\midrule
103 & $(1.35 \pm 0.04)\times10^{-7}$ & 1,025$\times$ & \textbf{7,407$\times$} \\
111 & $(1.82 \pm 0.08)\times10^{-8}$ & 1,691$\times$ & \textbf{13,671$\times$}\\
115 & $(6.2 \pm 0.6)\times10^{-9}$   & 2,309$\times$ & \textbf{22,888$\times$}\\
\bottomrule
\end{tabular}
\label{tab:rl_efficiency}
\end{table}

\section{Conclusions}
In this work, we have reformulated the Hadamard maximal determinant problem in orders $n \equiv 3 \pmod 4$ as a reinforcement learning based sequence reconstruction problem. Reconstruction of a candidate maxdet matrix then reduces to repeated sampling from the policy of a trained agent. The resulting constructions allow us to produce matrices in order $n=51$ and orders from $n=103$ to $n=119$ for $n\equiv 3\pmod 4$ greater than previously known largest determinants.

Although the size of the state space ranges from $10^{13}$ to $10^{33}$, decomposing the reconstruction problem into \textit{(a)} arithmetic feasibility and \textit{(b)} compatibility across different projections into cyclic groups, enables reinforcement learning to develop search strategies that generalize to previously unseen target joint autocorrelation functions. Furthermore, while classical methods, such as simulated annealing, restart search from scratch for each new target, a learned policy is able to leverage experience across samples, greatly reducing reconstruction time. This is an illustrative example of the potential benefits of RL for combinatorial optimization; amortization of inference costs over training leads to significantly lower per--instance costs relative to the best--performing classical methods.

\clearpage

\section*{Acknowledgements}
The project was sponsored by the Defense Advanced Research Projects Agency under cooperative agreement HR0011262E017, by the NSF AIMing grant 2522494, by the DRW Foundation, by a philanthropic gift from Les Kohn, and by a gift from Nebius Inc. X.H.\ was supported in part by an AMS--Simons Travel Grant and NSF AIMing grant 2617281. J.T. is also supported by a Cognia Junior fellowship. M.T. is also supported by the U.S. Department of Energy (Grant No. DE-SC0011632) and by the Walter Burke Institute for Theoretical Physics. Additionally, this work was supported with Cloud TPUs from Google's TPU Research Cloud (TRC), with GPUs from the NVIDIA Academic Grant Program, by cloud computing resources provided by Nebius through the Research Program of Nebius Academy, by Advanced Micro Devices, Inc. under the AMD University Program’s AI $\&$ HPC Cluster, and by computing resources from the Computer Lab at the University of Cambridge. The content of the information does not necessarily reflect the position or the policy of the Government, and no official endorsement should be inferred.

\bibliography{iclr2027_conference}
\bibliographystyle{iclr2027_conference}

\appendix

\section{Further results}\label{app:further_results}
\subsection{New record matrices}
We exhibit new record--determinant matrices obtained by reinforcement learning in orders $n=103,111,115$ in Table \ref{tab:record_matrices}. All matrices are of the bordered bicirculant form defined by the \textit{uncorrelated} construction \eqref{eq:ansatz}, so we elide the border, and only report the circulant component of the matrices $A_n = \text{circ}(a_n)$ and $B_n = \text{circ}(b_n)$ for each order $n$.

\begin{table}[t]
\centering
\small
\caption{The record-achieving matrices at orders $n = 2m + 1 = 103, 111, 115$. $\sigma_{\cdot}$ denotes the row sum, and
$\vert\det X\vert/D_{\mathrm{ref}}$ is the exact determinant over the best previously known determinant at that order \citep{orrick2012maxdet}. Determinants are computed using exact--integer arithmetic.}
\label{tab:record_matrices}
\def\pmseq#1{\texttt{\pmseqchar#1\pmseqend}}%
\def\pmseqend{\pmseqend}%
\def\pmseqchar#1{\ifx#1\pmseqend\else\ifx#1.\hspace{0.45em}\else#1\kern0pt\fi\expandafter\pmseqchar\fi}%
\smallskip
$\displaystyle X(A,B) = \begin{pmatrix} 1 & j^T & -j^T \\ j & A & B \\ -j & B^T & -A^T \end{pmatrix}, \qquad
A = \mathsf{Circ}(a), \quad B = \mathsf{Circ}(b), \qquad \mathsf{Circ}(x)_{ij} = x_{(j-i) \bmod m}$
\par\medskip
\begin{tabular}{@{}clr@{}}
\toprule
& First row & $\sigma_{\cdot}$ \\
\midrule
\multicolumn{3}{@{}l}{$n = 103$, \ $m = 51$, \ $\vert\det X\vert/D_{\mathrm{ref}} = 1.0000004$} \\
$a_{103}$ & \pmseq{--+++++++-.-+---+-++-.--++-+-+--.+-+-++--+-.++--+-+++-.-} & $1$ \\
$b_{103}$ & \pmseq{++--+--+-+.--+++-+-+-.--+-++--++.--++----++.++++++++-+.+} & $7$ \\
\addlinespace
\multicolumn{3}{@{}l}{$n = 111$, \ $m = 55$, \ $\vert\det X\vert/D_{\mathrm{ref}} = 1.0126294$} \\
$a_{111}$ & \pmseq{++++--++-+.---+-++-+-.+++++--++-.-+---++--+.-+---++++-.-+---} & $1$ \\
$b_{111}$ & \pmseq{+++-+--+--.++--+++-+-.---++-+-+-.+-++-+----.-++++-++++.-++++} & $7$ \\
\addlinespace
\multicolumn{3}{@{}l}{$n = 115$, \ $m = 57$, \ $\vert\det X\vert/D_{\mathrm{ref}} = 1.0182876$} \\
$a_{115}$ & \pmseq{++++----+-.-++-+-+--+.-+--+++---.+-+++-+--+.--++--+-++.+-+++--} & $1$ \\
$b_{115}$ & \pmseq{++----+-+-.-++++-+-+-.--++-++--+.----++-++-.--+++-+-++.+++++++} & $7$ \\
\bottomrule
\end{tabular}
\end{table}

\begin{table}[t]
\centering
\small
\caption{Record determinant sample statistics. `Distinct pairs' indicates pairs which are not related by a trivial joint symmetry to any other pair. `Samples needed' indicates the cumulative number of samples needed up to and including the first sampled batch which yielded a reconstruction of the highest--determinant target. The wallclock time is measured as the target's share of a batched run, i.e. (number of samples needed) $\times$ (average seconds per sample).}
\setlength{\tabcolsep}{6pt}
\begin{tabular}{rrrrrrrr}
\toprule
$n$ & \makecell[r]{$\det X/$\\record} & \makecell[r]{Recon-\\structions} & \makecell[r]{Distinct\\pairs} & \makecell[r]{Samples needed} & \makecell[r]{Samples\\drawn} & \makecell[r]{Time to\\first (s)} & \makecell[r]{Total\\time (s)} \\
\midrule
103 & 1.0000004 & 61 & 8 & $\leq 17{,}408$ & 148,480 & 1.2 & 10.4 \\
111 & 1.0126294 & 7 & 3 & $\leq 74{,}752$ & 271,360 & 18.9 & 66.7 \\
115 & 1.0182876 & 1 & 1 & $\leq 403{,}456$ & 534,528 & 44.0 & 56.9 \\
\bottomrule
\end{tabular}
\label{tab:record_samples}
\end{table}

\subsection{Sampling efficiency}
We exhibit the computational cost of finding the maximal determinant as predicted by the joint periodic autocorrelation functions in Table \ref{tab:record_samples}. For the RL approaches this reduces to sampling reconstructed sequences from a trained policy, conditioned on the target. Here we see that RL--based approaches allow us to amortize the cost of evaluation over training, and are significantly more efficient at inference---time than simulated annealing.

\section{Arithmetic--free baseline}\label{app:baseline_arch}

\subsection{Architecture}
Our architecture for the arithmetic--free baseline is a small autoregressive transformer decoder $D_{\varphi}$, which, in addition to an embedding of the previous action $a_{t-1}$, accepts as input the Fourier transform $\hat{p}_k$ of the target sequence, as well as an embedding of the last position committed $c_{t-1}$ and the remaining gaps to emit, $k_A - 1 -t$. The decoder consists of 3 attention blocks, $d = 192$, and 6 heads of width 32. $D_{\varphi}$ outputs logits corresponding to a sequence of valid gaps between consecutive $+1$ elements of $A$; $D_{\varphi}: \mathbb{R}^{(h+2) \times d} \rightarrow \mathbb{R}^T$, where $T = k_A - 1$ is the trajectory length. The policy is an MLP head over the output of the decoder. The dataset for each order consists of all possible enumerated periodic autocorrelation targets $p(t)$ for that order, of size $O(10^5)$. We order the dataset by determinant and train on the bottom 95\% of sequences, retaining the top 5\% of sequences for evaluation. The policy is trained using PPO without an entropy term in the objective, with $\alpha=0.03$ set after a sweep at $n=103$. We set $\gamma=1$ and the GAE lookahead $\lambda=0.95$ without tuning. The critic is an MLP head over the objective and shares the same decoder trunk as the policy head. The architecture is held identical across all orders apart from the necessary changes to input dimensions and embedding tables.

\subsection{Dataset}
To generate the dataset used for the arithmetic--free method, note the orders are sufficiently small for exhaustive enumeration of  all possible autocorrelation functions\footnote{$p(t)$ is a sequence of length $m=(n-1)/2$ with reflection symmetry along the midpoint, so one only needs to enumerate all $(m-1)/2$--sequences of binary strings, compute their orbit under the permutation $p \mapsto p \circ \lambda$ for some multiplicative index shift $\lambda \in \mathbb{N}$, and record the canonical representative of this orbit as the lexicographic minimum.}. We restrict ourselves to targets which yield a determinant of $80\%$ or higher of the previous record for each order, which yields an $O(10^4)$ dataset $\mathcal{D}_n = \{p_n(t)\}$ of autocorrelation functions at all orders $n \equiv 3\pmod 4$ we consider, which form a superset of the feasible autocorrelation functions admitting the decomposition \eqref{eq:paf_constraint}.

\subsection{Sampling efficiency}
We exhibit further results showcasing the efficiency of a trained policy over a random baseline in Table \ref{tab:validation_samples}. Note we have recast the problem into a form that is nominally insensitive to the determinant of the bicirculant matrix $X$ constructed from the pair $(A,B)$, reflected in the pass@1 rate being approximately constant for the top 5\% and top 0.5\% of determinant--maximizing targets. 

\begin{table}[t]
\centering
\small
\caption{Pair $(A,B)$ reconstruction statistics on the validation set $V$ (top 5\% of targets measured by determinant) and on its top end (top 0.5\%) for the arithmetic--free baseline . `Fraction reconstructed' counts targets with at least one certified pair over all samples drawn. `Samples/target' is the largest number drawn on any one target. pass@$k$ is the mean over targets drawn from the respective sets of the unbiased estimator $\hat{\rho}_k$ \eqref{eq:unbiased_est_rho_k}, while pass@1 is the mean of $n^+_i/n_i$ over $V$.}
\label{tab:validation_samples}
\setlength{\tabcolsep}{1.5pt}
\begin{tabular}{rrrrrrrr}
\toprule
$n$ & Reconstructed & Samples/target & \makecell[r]{pass@1\\[-2pt]{\scriptsize$\times10^{-5}$}} & \makecell[r]{pass@256} & \makecell[r]{pass@1024} & \makecell[r]{pass@$2^{12}$} & \makecell[r]{pass@$2^{16}$} \\
 &  & \emph{\scriptsize max.\ MC error} & {\scriptsize$\pm5.4\%$} & {\scriptsize$\pm5.4\%$} & {\scriptsize$\pm5.3\%$} & {\scriptsize$\pm4.9\%$} & {\scriptsize$\pm2.5\%$} \\
\midrule
\multicolumn{8}{l}{\emph{Validation set (top 5\%)}} \\
103 & 560/560 & $\le$345k & 13.8 & 0.0346 & 0.130 & 0.414 & 0.994 \\
111 & 744/744 & $\le$402k & 3.07 & 0.00783 & 0.0309 & 0.118 & 0.849 \\
115 & 1055/1057 & $\le$2.6M & 1.43 & 0.00366 & 0.0145 & 0.0565 & 0.572 \\
\midrule
\multicolumn{8}{l}{\emph{Top 0.5\%}} \\
103 & 56/56 & \multirow{3}{*}{(same as above)} & 15.6 & 0.0389 & 0.145 & 0.447 & $>$0.999 \\
111 & 74/74 &  & 3.13 & 0.00797 & 0.0315 & 0.120 & 0.864 \\
115 & 106/106 &  & 1.43 & 0.00366 & 0.0146 & 0.0566 & 0.573 \\
\bottomrule
\end{tabular}
\end{table}

\section{Further results on PAFs}\label{app:paf}

\begin{prop} Let $N$ be the number of times $p(s) = +2$ for $s\neq 0$, then:
\begin{gather}
    |a|^2 + |b|^2 = 4N+2\,.
\end{gather}
\end{prop}
\begin{proof} The proof of the original statement can be found in \citep{djokovic2011new}. Let $N$ be the number of times $p(s) = 2$ for $s\neq 0$, then:
\begin{gather}\label{eq:pafSum}
    \sum_{s=0}^{m-1} p(s) = 2m + 2N - 2(m - 1 - N) = 4N + 2\,.
\end{gather}
Using the definition of the periodic autocorrelation functions \eqref{eq:paf}, we see that:
\begin{gather}\label{eq:pafSum2}
    \sum_{s=0}^{m-1} \mathrm{PAF}_a(s) = \sum_{s=0}^{m-1}a_s^2 + \sum_{j=0}^{m-1}\sum_{s=1}^{m-1}a_j a_{j+s} = \sum_{s=0}^{m-1}a_s^2 + 2\sum_{j<k}a_j a_k = \left(\sum_{j=0}^{m-1}a_j\right)^2
\end{gather}
Let $|a|, \, |b|$ denote the row sum of $a$, $b$, respectively. Combining \eqref{eq:pafSum} and \eqref{eq:pafSum2} yields $|a|^2 + |b|^2 = 4N + 2$.
\end{proof}

\subsection{Correlated $X^{\mathrm{corr}.}$}\label{app:corr}
 We follow the same exact procedure for $X = X^{\mathrm{corr}.}(A,B)$. Except, the Gram matrix $G$ of $X$ is given by:
\begin{gather}
    G = XX^T = \begin{pmatrix}
        n & \gamma^T & \gamma^T \\
        \gamma & C & j\otimes j - K\\
        \gamma & j\otimes j-K^T & C
    \end{pmatrix}\,,
\end{gather}
where $\gamma = (1+|a|-|b|)j$ and $C$ is same as before. The main structural difference is the off-diagonal contribution coming from $K = A(B+B^T)$. Following Ehlich Gram constraint, we impose $K_{ij} = 2$ condition on $A$ and $B$. While it is impossible to rewrite this condition purely in terms of PAF-sum targets $p(s)$, we may still write it in terms of the constraints on the cross-correlations of $a$ and $b$. In particular, by taking the Fourier transform of $K = 2j\otimes j$ it is easy to show that the relevant constraints are:
\begin{gather}
    \text{Zero mode:}~|a||b| = m,\quad \text{Higher modes:}~~\widehat{a}_k\mathrm{Re}(\widehat{b}_{k}) = 0\,,
\end{gather}
where latter may be rewritten purely in terms of constraints on the convolution of sequences: $a* c = 2|a||b|/m = 2$, where $c_t := b_{t} + b_{-t\pmod m}$. Assuming such sequences exist, we calculate the determinant in the same manner as before: the Fourier modes of $X$ are given by:
\begin{gather}
    \widehat{X}_0 = \begin{pmatrix}
        1 & \sqrt{m} & -\sqrt{m} \\
        \sqrt{m} & |a| & |b| \\
        \sqrt{m} & -|b| & -|a|
    \end{pmatrix},\quad \widehat{X}_{t>0} = \begin{pmatrix}\widehat{a}_t & \widehat{b}_t \\
    -\widehat{b}_t & -\overline{\widehat{a}}_t\end{pmatrix}\,.
\end{gather}
While $\det \widehat{X}_0$ similarly only depends on the segment sums $|a|$ and $|b|$, the form of the determinant of the higher modes $\det\widehat{X}_{t>0}$ is given by:
\begin{gather}
    |\det\widehat{X}_{t>0}|^2 = \left|\widehat{b}_t^2-|\widehat{a}_t|^2\right|^2 = (|\widehat{a}_t|^2 + |\widehat{b}_t|^2)^2 -4 (|\widehat{a}_t|\mathrm{Re}(\widehat{b}_t))^2 = (|\widehat{a}_t|^2 + |\widehat{b}_t|^2)^2 
\end{gather}
where the third equality holds true due to assumption $K = 2j\otimes j$. Using the same expansions of higher mode determinants as before, we obtain the same exact expression for the determinant of the higher modes of $X$, while zero mode can similarly be determined from the segment sums. The exhaustive list of all admissible targets for $n=51$ parametrized using $X^{\mathrm{corr}.}$ is shown in Fig.~\ref{fig:paftargets}.

\section{Numerical and RL approaches}

\begin{table}[h!]
    \centering
    \begin{tabular}{@{}ccc@{}}
         \toprule
         Order $n$ & Segment sums $(|a|,|b|)$ & $\log_{10}\left|\mathcal{S}(|a|,|b|,(n-1)/2)\right|$ \\
         \midrule
         $51$ & $(5,5)$ & 13.03 \\
         $103$ & $(1,7)$ & 28.59\\
         $107$ & $(1,7)$ & 29.78\\
         $111$ & $(1,7)$ & 30.98\\
         $115$ & $(1,7)$ & 32.17\\
         $119$ & $(1,7)$ & 33.36\\
         \bottomrule
    \end{tabular}
    \caption{Size of the state space $\mathcal{S}(|a|,|b|,(n-1)/2)$ for orders $n$ where we search for $n=51$ and $103\leq n\leq 119$ with $n\equiv 3\pmod 4$.}
    \label{tab:stateSpaceScaling}
\end{table}

\subsection{Order $51$ using $X^{\mathrm{corr}.}(A,B)$}
Simulated annealing with energy functional $E$ given by:
\begin{gather}
    E(a,b) := H(a,b) + \sum_{t=0}^{m-1}\bigg|(a*c)_t - 2\bigg|^2
\end{gather}
Where $c_t = b_t + b_{-t\pmod m}$ and $H(a,b)$ is defined in \eqref{eq:energy}. We reconstructed all 10 record--breaking PAF--sum targets at order 51. The reconstructions are summarized in Table~\ref{table:recon51}.

\begin{table}[h!]
    \centering
        \small
    \setlength{\tabcolsep}{8pt}
    \renewcommand{\arraystretch}{1.25}
    \resizebox{\linewidth}{!}{%
        \begin{tabular}{@{}lll@{}}
            \toprule
            \multicolumn{1}{c}{Sequence $a$} &
            \multicolumn{1}{c}{Sequence $b$} &
            \multicolumn{1}{c}{PAF-sum target} \\
            \midrule
            \signseq{--+-+-+---++++++-++-++-++} &
            \signseq{+++++++-+++--+++--+-+----} &
            \texttt{+2 +2 +2 +2 +2 +2 -2 -2 -2 -2 -2 -2} \\

            \signseq{-+++-+++-++-+-----++++-++} &
            \signseq{+----++-+-+++--++-+-+++++} &
            \texttt{+2 +2 +2 -2 -2 -2 -2 -2 -2 +2 +2 +2} \\

            \signseq{++-+-++-----+-++++++--+++} &
            \signseq{++++-++--+++-+-+-++-++---} &
            \texttt{+2 +2 -2 -2 -2 -2 +2 +2 +2 +2 -2 -2} \\

            \signseq{++++-++++--+--++-+++----+} &
            \signseq{++-+-+---++--+++-+++++-+-} &
            \texttt{+2 -2 -2 +2 +2 -2 -2 +2 +2 -2 -2 +2} \\

            \signseq{+++----++++++-++--++-+-+-} &
            \signseq{+++--+-+--+--+++++-++++--} &
            \texttt{+2 -2 -2 -2 +2 +2 +2 -2 -2 -2 +2 +2} \\

            \signseq{+-++-+++-+-+-----+++++-++} &
            \signseq{+++-+++--++++--+-++-+-+--} &
            \texttt{-2 +2 -2 +2 +2 -2 +2 -2 +2 -2 +2 -2} \\

            \signseq{+++-+-+++--++-++--+-+--++} &
            \signseq{+----+-+-++--+++-+-++++++} &
            \texttt{-2 +2 -2 +2 -2 +2 -2 +2 -2 +2 -2 +2} \\

            \signseq{+++++-++--+++-+---+-+--++} &
            \signseq{++-+-+++--+--+++++--++-+-} &
            \texttt{-2 -2 +2 +2 -2 -2 +2 +2 -2 +2 +2 -2} \\

            \signseq{+----++-++--+-+-+++-+++++} &
            \signseq{+---++--+-++-+-++-+++-+++} &
            \texttt{-2 -2 +2 -2 +2 +2 -2 +2 +2 -2 +2 -2} \\

            \signseq{+-++++--+++++-+-+++--+---} &
            \signseq{++---+-++-++-+-++--+++++-} &
            \texttt{-2 -2 +2 -2 -2 +2 +2 -2 +2 +2 -2 +2} \\
            \bottomrule
        \end{tabular}%
    }
    \caption{Exact reconstructions for all $10$ maximal determinant PAF-sum targets for order $n=51$.}\label{table:recon51}
\end{table}

We limited experimentation to simulated annealing in this order as it was able to produce solutions almost instantaneously.

\subsection{Orders $107$ and $119$ using simulated annealing}
Since $m=(n-1)/2$ for both orders $107$ and $119$ is prime, the lifting setup described in Sec~\ref{ssec:oracle} does not apply. However, simulated annealing on the energy functional $E(a,b) := H(a,b)$ \eqref{eq:energy} yields the maximal determinant matrices after running for few days. We report reconstruction of a maximal determinant PAF-sum target for each of the orders $107$ and $119$ in Table~\ref{table:recon107119}.

\begin{table}[h!]
    \centering
        \small
    \setlength{\tabcolsep}{8pt}
    \renewcommand{\arraystretch}{1.25}
    \resizebox{\linewidth}{!}{%
    \begin{tabular}{@{}clll@{}}
        \toprule
        Order $n$ &
        \multicolumn{2}{c}{Sequences} &
        \multicolumn{1}{c}{PAF-sum target} \\
        \midrule

        \multirow{2}{*}{$107$}
        & $a$ &
        \signseq{-+-+-++++-+--++-+-++++--++++-++--------+-++-++-----++}
        & \multirow{2}{*}{%
            \shortstack[l]{%
                \texttt{-2 +2 -2 +2 -2 -2 -2 -2 -2 -2 -2 -2 -2}\\
                \texttt{-2 +2 -2 +2 -2 +2 -2 +2 -2 -2 -2 -2 -2}%
            }} \\
        & $b$ &
        \signseq{++-+-+---++++-+-++++-++-+++---+---++--++-++-+++-+---+}
        & \\

        \midrule

        \multirow{2}{*}{$119$}
        & $a$ &
        \signseq{-----+---+-++-+-+-+--+--+++++++-+-++-++----++-+--++--++-+++}
        & \multirow{2}{*}{%
            \shortstack[l]{%
                \texttt{-2 -2 -2 +2 -2 -2 -2 +2 -2 -2 -2 -2 +2 -2 -2}\\
                \texttt{-2 +2 -2 -2 -2 +2 -2 -2 -2 +2 -2 -2 -2 -2}%
            }} \\
        & $b$ &
        \signseq{++++++++---+-+++--+++--++-+-++++---+--++-+-+----++-+---+++-}
        & \\

        \bottomrule
    \end{tabular}%
    }
    \caption{Exact reconstructions of a maximal determinant PAF-sum target for each order $n\in\{107,119\}$.}\label{table:recon107119}
\end{table}

\subsection{Order $103$ using RL on full state space $\mathcal{S}(|a|,|b|,m)$}\label{sec:jointRL}
\begin{figure}
    \centering
    \includegraphics[width=0.8\linewidth, trim = {0cm, 0.5cm, 0cm, 0cm}]{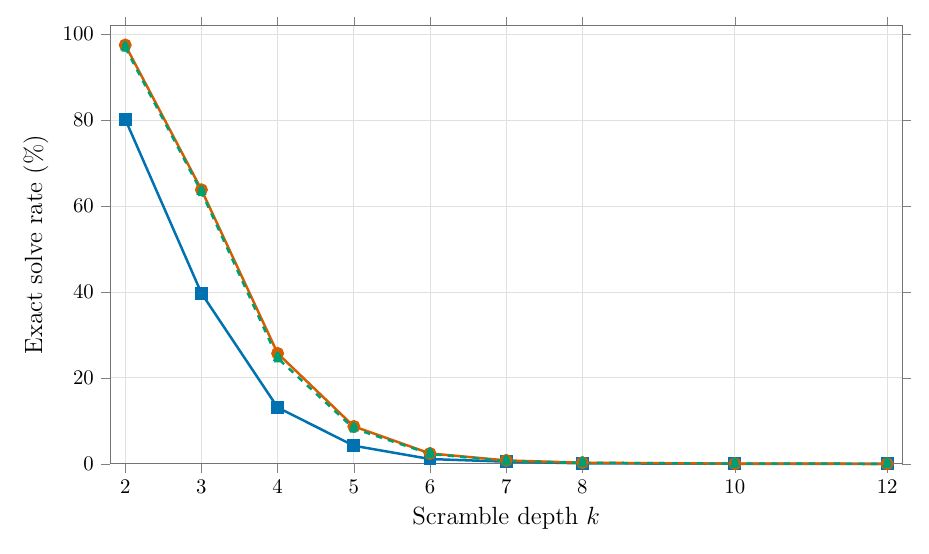}
    \caption{Exact reconstruction rate of non-arithmetic energy-based RL on the full state space $\mathcal{S}(|a|,|b|,m)$ for held-out order-103 targets after $k$ random steps, with a budget of $2k$ moves.
    Energy-guided search (blue), imitation learning (orange), and GRPO fine-tuning (green, dashed), all decline sharply as scramble depth increases}
    \label{fig:reconstruction rates}
\end{figure}

From an exhaustive enumeration of PAF-sum targets in $n\equiv 3\pmod 4$ (illustrated in Fig.~\ref{fig:paftargets}), maximal--determinant candidates in dimensions from $103$ to $119$ whose determinant exceeds the previous record \citep{orrick2012maxdet} arise from the $X(A,B)$ family. For all of these matrices the segment sums attaining the maximal determinant are $(1,7)$. Since these matrices come from $X(A,B)$, in order for a pair of sequences $(a,b)$ to realize those targets, they must solve only one constraint:
\begin{gather}
    H(a,b) = 0\,
\end{gather}
where $H(a,b)$ is the PAF-sum objective defined in \eqref{eq:energy}. In this section, we focus specifically on order $103$ and describe the training algorithm of an agent which takes an arbitrary PAF-sum target as an input and produces its reconstruction. Motivated by the results of \citep{ramrakhya2023pirlnav,agostinelli2019solving, huangposition}, we first sample a small-scale PAF-sum reconstruction dataset (Alg.~\ref{alg:samplePretrainingData}), pretrain the agent on it using behavior cloning (Alg.~\ref{alg:pretrainAgent}), and then finetune it through a curriculum of increasing scramble depth (Alg.~\ref{alg:curriculumFinetuning}) to allow it to generalize to reconstruct unseen PAF-sum targets. Conceptually, the illustration of this transfer and generalization to interesting PAF-sum targets is shown in Fig.~\ref{fig:pafRL}. We pick $\mathcal{D}_{\mathrm{val}}$ to consist of top $20\%$ of all enumerated PAF-sum targets in Fig.~\ref{fig:paftargets}. At inference, the trained agent attempts to reconstruct these targets from uniformly sampled initial states (Alg.~\ref{alg:inferenceGeneration}). The PAF-sum targets $\mathcal{D}_{\mathrm{PAF}}$ used for pretraining and finetuning do not overlap with $\mathcal{D}_{\mathrm{val}}$, moreover, it is highly unlikely that it will contain reconstruction of the maximal PAF-sum target. This is mainly because constructing an admissible PAF-sum target is a significantly easier problem than reconstructing a specific PAF-sum target.
Figure \ref{fig:reconstruction rates} compares reconstruction rates for energy-guided, imitation learning and GRPO fine tuning for a $k$ random swap scramble.
We observe that at scramble depth 10-12 the imitation learning policy is unable to reconstruct the PAF-sum making it unfeasible for constructing full solutions.

\subsection{Pretraining and Curriculum algorithms}\label{app:pt_curr_alg}
\begin{algorithm}[h!]
\caption{Sampling pretraining PAF data}\label{alg:samplePretrainingData}
\begin{algorithmic}[1]
    \Require Order $n\equiv 3\pmod 4$, $m=(n-1)/2$. Pretraining dataset size $N>0$.
    \Statex \textbf{I: Sampling pretraining data}
    \State $\mathcal{D}_{\mathrm{PAF}}\gets \emptyset$
    \While{$|\mathcal{D}_{\mathrm{PAF}}| < N$}
        \State $(a, b)\gets \Call{SimulatedAnnealing}{m}$
        \State $p \gets \mathrm{PAF}_{a} + \mathrm{PAF}_{b}$
        \If{$\forall s\neq 0\colon p(s)\in \{-2,+2\}$}
            \State $\mathcal{D}_{\mathrm{PAF}} \gets \mathcal{D}_{\mathrm{PAF}}\cup\{(a, b)\}$
        \EndIf
    \EndWhile
\end{algorithmic}
\end{algorithm}

\begin{algorithm}[h!]
\caption{Pretraining the agent}\label{alg:pretrainAgent}
\begin{algorithmic}[1]
\Require Order $n\equiv 3\pmod 4$. Pretraining dataset $\mathcal{D}_{\mathrm{PAF}}$. Max scrambles $K_{\mathrm{max}}>0$. Swap-action environment $\mathsf{env}$.
\State $\mathcal{D}_{\mathrm{traj}}\gets\emptyset$
\For{$(a, b, p)\in \mathcal{D}_{\mathrm{PAF}}$}
\State $s\gets (a,b)$
\For{$k=1,\ldots,K_{\mathrm{max}}$}
    \State $(s_0 = s, s_1,\dots, s_k) \gets \Call{RandomScramble}{s,k,\mathsf{env}}$
    \State $\tau \gets \Call{ShortestPath}{s_0, s_k}$
    \State $\mathcal{D}_{\mathrm{traj}}\gets \mathcal{D}_{\mathrm{traj}}\cup \{(\tau, p)\}$
\EndFor
\EndFor
\State $\theta \gets \Call{BC}{\pi_\theta, \mathcal{D}_{\mathrm{traj}}}$\Comment{Behavior Cloning}
\end{algorithmic}
\end{algorithm}

\begin{algorithm}[h!]
\caption{Curriculum learning over increasing scramble depth}\label{alg:curriculumFinetuning}
\begin{algorithmic}[1]
\Require Order $n\equiv 3\pmod 4$. Pretraining dataset $\mathcal{D}_{\mathrm{PAF}}$. Swap-action environment $\mathsf{env}$. Curriculum scheduler $k_t$.
\For{$t=1,\ldots, T_\mathrm{max}$}
    \State $(p,s)\gets \mathcal{D}_{\mathrm{PAF}}$
    \State $\epsilon \sim \mathcal{N}(0,1)$
    \If{$\epsilon > \epsilon_{\mathrm{th}}$}
    \State $(\dots, s_k) \gets \Call{RandomScramble}{s,k_t,\mathsf{env}}$
    \Else
        \State $s_k \sim \mathcal{U}(\mathsf{env})$\Comment{Random uniform from state space.}
    \EndIf
    \State $\tau \gets \Call{Rollout}{\pi_\theta, s_k, p, \mathsf{env}}$
    \State $\theta \gets \Call{RLOptimize}{\theta,\tau, p}$
\EndFor
\end{algorithmic}
\end{algorithm}

\begin{algorithm}[h!]
\caption{Inference / generation}\label{alg:inferenceGeneration}
\begin{algorithmic}[1]
    \Require Order $n\equiv 3\pmod 4$. Max scrambles $K_{\mathrm{max}}>0$. Swap-action environment $\mathsf{env}$. Trained policy $\pi_\theta$.
    \State {\bf Inference / generation}
    \State $S\gets \emptyset$
    \For{$p\in \mathcal{D}_{\mathrm{val}}$}
        \State $s_0\gets \mathcal{U}(\mathsf{env})$
        \State $s^\star \gets \Call{Rollout}{\pi_\theta,s_0,p,\mathsf{env}}$
        \If{$H(s^\star) = 0$}
            \State $S\gets S\cup \{(s^\star, p)\} $
        \EndIf
    \EndFor
    \State \Return $\pi_\theta, S$
\end{algorithmic}
\end{algorithm}

\clearpage
\section{Arithmetic Obstructions}\label{app:arithmetic}
The key identification is that we may regard sequences $\{a_j\}_{j=0}^{m-1}$ as elements of $\mathbb{Q}[C_m]$ where $C_m$ is a cyclic group of order $m$. This therefore allows us to interpret PAF functions of $a$ as products:
\begin{gather}
    aa^* = \sum_{j=0}^{m-1}\mathrm{PAF}_a(j)z^j\in\mathbb{Q}[z]/(z^m - 1)\cong \mathbb{Q}[C_m]
\end{gather}
Therefore, the problem of PAF reconstruction of some profile $r$ is essentially a problem of searching for $a\in\mathbb{Q}[C_m]$ such that $aa^* = r$, where $r$ here is implicitly identified as an element of $\mathbb{Q}[C_m]$. While this reformulation doesn't necessarily simplify the problem, we may look for local obstructions by projecting $a$ to the cyclotomic subfields $\mathbb{Q}(\zeta_d)$ of the decomposition:
\begin{gather}
    \mathbb{Q}[C_m] = \prod_{d\mid m} \mathbb{Q}(\zeta_d)\,,
\end{gather}

\begin{thm} Let $d\mid m$ for some $d>2$ and $K_d^+$ be the maximal real subfield of $K_d := \mathbb{Q}(\zeta_d)$, then $b_db_d^* = r_d\in (K_d^+)^\times$ has a solution $b_d\in K_d$ iff:
\begin{gather}\label{eq:cohomology}
    [r_d] = 0~\text{in}~H^2(K_d/K_d^+)\,,
\end{gather}
where $H^2(K_d/K_d^+)$ is the Galois cohomology $H^2(\mathrm{Gal}(K_d/K_d^+), K_d^\times)$.
\end{thm}
\begin{proof} Since $H^2(K_d/K_d^+) \simeq (K_d^+)^\times / N_{K_d/K_d^+}(K_d^\times )$ \citep{serre1979local}. Let $c$ be the generator of $\mathrm{Gal}(K_d/K_d^+)\simeq \mathbb{Z}/(2)$, then we have:
\begin{gather}
    [r_d] = 0~\Longleftrightarrow~ r_d = N_{K_d/K_d^+}(b)
\end{gather}
for some $b_d\in K_d$, thus $N_{K_d/K_d^+}(b_d) = c(b_d)b_d = b_db_d^* = r_d$.
\end{proof}
The case $r_d = 0$ can be handled separately, as it has the unique solution $b_d=0$.

We may infer existence of global obstruction $[r_d]\neq 0$ by doing finitely many local tests. In particular, we have:
\begin{thm}\label{thm:obstructions} Let $\Delta_d := (\zeta_d-\zeta_d^{-1})^2$. Let $r_d$ be such that it passes the Archimedean test: $\sigma(r_d) > 0$ for all embeddings $\sigma\colon K_d^+\rightarrow \mathbb{R}$. Then, we have the following equivalence:
\begin{gather}\label{eq:hilbert}
    \forall 0\neq v\in\mathrm{Spec}~\mathcal{O}_{K_d^+}\colon (r_d, \Delta_d)_v = 1 ~\Longleftrightarrow~ [r_d]=0~\text{in}~H^2(K_d/K_d^+)\,,
\end{gather}
where $(r_d,\Delta_d)_v$ is a Hilbert symbol at a place $v$.
\end{thm}
\begin{proof} The element $\zeta_d\in K_d$ has a minimal polynomial over $K_d^+$ given by $f_d = X^2 - s_dX + 1=0$, where $s_d = \zeta_d+\zeta_d^{-1}$. The discriminant of $f_d$ is $\Delta_d = s_d^2 - 4 = (\zeta_d-\zeta_d^{-1})^2$, thus $K_d = K_d^+(\sqrt{\Delta_d})$. From \citep{Serre:1993}, the condition $(r_d,\Delta_d)_v = 1$ is equivalent to $r_d \in N_{(K_d)_v^+(\sqrt{\Delta})/(K_d)_v^+}((K_d)_v^+(\sqrt{\Delta_d})^\times)$, i.e. $r_d$ is a local norm. Since $(r_d,\Delta_d)_v = 1$ at all finite places $v$ and $\sigma(r_d)>0$ for all $\sigma\colon K_d^+\rightarrow \mathbb{R}$, by Hasse norm theorem, $r_d$ must also be global norm. Converse is obvious.
\end{proof}

We use \texttt{cypari2} \citep{cypari2, PARI2} to efficiently calculate Hilbert symbols \eqref{eq:hilbert} at all relevant places $v\in\mathrm{Spec}~\mathcal{O}_{K_d^+}$.

\end{document}